\documentclass[preprint,11pt]{elsarticle}

\newtheorem{theorem}{Theorem}[section]
\newtheorem{lemma}{Lemma}[section]

\newtheorem{proposition}{Proposition}[section]

\newcommand{\ignore}[1]{}{}
\newproof{pf}{Proof}

\usepackage{amssymb}
\usepackage{amsmath}
\usepackage{amsfonts,graphicx,color}

\usepackage{color}
\usepackage{soul}
\usepackage[dvipsnames]{xcolor}

\usepackage[colorlinks=true, linkcolor=blue, citecolor=blue]{hyperref}

\newcommand\beq{\begin{equation}}
	\newcommand\eeq{\end{equation}}

\usepackage{ifthen}
\usepackage{xkeyval}
\usepackage{todonotes}
\makeatletter

\newcommand{\Rmnum}[1]{\expandafter\@slowromancap\romannumeral #1@}

\begin{document}
	\let\today\relax
	
\begin{frontmatter}
		\date{\quad}
		
		\title{Wasserstein-$1$ distance and nonuniform Berry-Esseen bound  for a supercritical branching process in a random environment}
	\author[cor1]{Hao Wu}
	\author[cor2]{Xiequan Fan}
	\author[cor3]{Zhiqiang Gao}
	\author[cor4]{Yinna Ye\corref{c1}}
	\ead{yinna.ye@xjtlu.edu.cn}
		
		\cortext[c1]{Corresponding author}
		\address[cor1]{Center for Applied Mathematics, Tianjin University, Tianjin 300072, P. R. China}
		\address[cor2]{School of Mathematics and Statistics, Northeastern University at Qinhuangdao, Qinhuangdao 066004, P. R. China}
		\address[cor3]{Laboratory of Mathematics and Complex Systems, School of Mathematical Sciences, Beijing Normal University, Beijing $100875 $, P. R. China}
		\address[cor4]{Department of Applied Mathematics, School of Mathematics and Physics, Xi'an Jiaotong-Liverpool University, Suzhou $215123$, P. R. China}
				
		\begin{abstract}
		Let $ (Z_{n})_{n\geq 0} $ be a supercritical branching process in an independent and identically distributed random environment.  We establish an optimal convergence rate in the Wasserstein-$1$ distance  for the process $ (Z_{n})_{n\geq 0} $, which completes a result of Grama et al. [Stochastic Process. Appl., \textbf{127}(4), 1255-1281, 2017]. Moreover, an exponential nonuniform Berry-Esseen bound is also given. At last, some applications  of the main results to the confidence interval estimation for the criticality parameter and the population size $Z_n$ are discussed.
		\end{abstract}
		
		\begin{keyword} Branching processes; Random environment; Wasserstein-$1$ distance;  Nonuniform Berry-Esseen bounds
			\vspace{0.3cm}
			\MSC 60J80; 60K37; 60F05; 62E20
		\end{keyword}
  \end{frontmatter}

\section{Introduction}
\setcounter{equation}{0}

Branching process in a random environment (BPRE) initially introduced by Smith and Wilkinson \cite{smith1969branching} is a generalization of the Galton-Watson process.  Denote by $\xi=\left(\xi_{0}, \xi_{1}, \ldots\right)$ a sequence of independent and identically  distributed (i.i.d.) random variables,  where $ \xi_{n} $ stands for the random environment in the $n$-th generation. Then the branching process $ (Z_{n})_{n\geq 0} $ in the random environment $ \xi $ is defined as follows:
$$
Z_{0}=1, \ Z_{n+1}=\sum_{i=1}^{Z_{n}} X_{n, i}, \ n \geq 0,
$$
where $X_{n, i}$ represents the number of offspring produced by the $i$-th particle in the $n$-th generation.  The distribution of $ X_{n,i} $ depending on the environment $ \xi_{n} $ is denoted by $p\left(\xi_{n}\right)=\left\{p_{k}\left(\xi_{n}\right)=\mathbb{P}(X_{n,i}=k\; |\xi_{n}): k \in \mathbb{N}\right\}$. Suppose that given $ \xi_{n} $, $ (X_{n, i})_{i\geq1}$ is a sequence of i.i.d. random variables; moreover, $ (X_{n, i})_{i\geq1}$ is independent of $(Z_{1},\ldots ,Z_{n}) $. Let $\left(\Gamma , \mathbb{P}_{\xi}\right)$ be the probability space under which the process is defined when the environment $\xi$ is given, where $\mathbb{P}_{\xi}$ is usually called the quenched law. The state space of the random environment $ \xi $ is denoted by $ \Theta $ and the total probability space can be regarded as the product space $\left(\Theta^{\mathbb{N}}\times\Gamma, \mathbb{P}\right),$ where $\mathbb{P}(\mbox{d} x, \mbox{d} \xi)=\mathbb{P}_{\xi}(\mbox{d} x) \tau(\mbox{d} \xi).$ That is, for any measurable positive function $g $ defined on $\Theta^{\mathbb{N}}\times\Gamma$, we have
$$\int g(x, \xi) \mathbb{P}(\mbox{d} x, \mbox{d} \xi)=\iint g(x, \xi) \mathbb{P}_{\xi}(\mbox{d} x) \tau(\mbox{d} \xi),$$
where $ \tau $ represents the distribution law of the random environment $ \xi $. We usually call $ \mathbb{P} $ the annealed law. And $\mathbb{P}_{\xi} $ can be regarded as the conditional probability of $\mathbb{P} $ given the environment $ \xi $.
The expectations with respect to $ \mathbb{P}_{\xi} $ and $ \mathbb{P} $ are denoted by $ \mathbb{E}_{\xi} $ and $\mathbb{E}$, respectively. For any environment $\xi$, integer $n\geq0$ and real number $p>0$, define
$$m_n^{(p)}=m_n^{(p)}(\xi)=\sum_{i=0}^{\infty}i^pp_i(\xi_n),\quad m_n=m_n(\xi)=m_n^{(1)}(\xi).$$
Then
$$m_0^{(p)}=\mathbb{E}_{\xi} Z_1^p\quad\text{and}\quad m_n=\mathbb{E}_{\xi}X_{n, i},\quad \ i\geq 1.$$
Consider the following random variables
$$\Pi_{0}=1,\ \Pi_{n}=\Pi_{n}(\xi)=\prod_{i=0}^{n-1} m_{i},\quad n\geq 1.$$
It is easy to see that
$$\mathbb{E}_{\xi}Z_{n+1}
=\mathbb{E}_{\xi}\Big[\sum_{i=1}^{Z_{n}} X_{n,i} \Big]\\
=\mathbb{E}_{\xi}\Big[\mathbb{E}_{\xi}\Big(\sum_{i=1}^{Z_{n}} X_{n,i} \Big|Z_{n}\Big)\Big]\\
=\mathbb{E}_{\xi}\Big[\sum_{i=1}^{Z_{n}} m_{n} \Big]\\
=m_{n}\mathbb{E}_{\xi}Z_{n}.$$
Then by recursion, we get
$\Pi_{n}= \mathbb{E}_{\xi}Z_{n}.$
Denote by $$ X=\log m_{0},\ \ \ \mu=\mathbb{E}X \ \ \textrm{and} \ \ \sigma^2=\mathbb{E}(X-\mu)^2. $$
The branching process $  (Z_{n})_{n\geq0} $ is called supercritical, critical or subcritical according
to $ \mu > 0 $, $ \mu = 0 $ or $  \mu < 0 $, respectively.   Hence $ \mu $ is known as the criticality parameter. Over all the paper, assume that 
\begin{equation}\label{H2}
	p_0\left(\xi_0\right)=0\text{ a.s.}
\end{equation}
The condition (\ref{H2}) means that each particle has at least one offspring, which implies $X\geq0$ a.s. and hence $Z_n\geq1$ a.s. for any $n\geq1$. 

Limit theorems for  BPRE have attracted a lot of interests since the seminal work of Smith and Wilkinson \cite{smith1969branching}.
For critical and subcritical BPRE, the study  mainly focuses on the survival probability and conditional limit theorems for the branching processes, see, for instance,  Vatutin \cite{Vatutin} and Afanasyev et al. \cite{Afanasyev,Boinghoff}. For the supercritical BPRE, a
number of researches focused on the central limit theorem, moderate and large deviations; see  B{\"o}inghoff and Kersting \cite{boinghoff2010upper}, Huang and Liu \cite{huang2012moments},
B{\"o}inghoff\ \cite{boinghoff2014limit}, Li et al. \cite{li2014asymptotic}, Grama et al. \cite{grama2017berry},
Wang and Liu\ \cite{wang2017limit}, Fan et al. \cite{fan2020uniform} and Gao \cite{G2021}.

The aim of this paper is to establish optimal convergence rates in the Wasserstein-$1$ distance for the branching process  in a random environment.  We first recall the definition of the Wasserstein-$1$ distance.  Let $\mathcal{L}(\mu, \nu)$ be a set
of probability laws on $\mathbb{R}^2$ with marginals $\mu$ and  $\nu$.  The Wasserstein-$1$ distance between $\mu$ and  $\nu$  is defined
as follows:
\begin{eqnarray*}
	W_1(\mu, \nu)= \inf  \bigg \{  \int |x-y|  \mathbb{P}(\mbox{d}x, \mbox{d}y)    :\ \mathbb{P} \in \mathcal{L}(\mu, \nu)      \bigg\} .
\end{eqnarray*}
In particular, if $\mu_X$ is  the distribution of $X$  and $\nu$ is the standard normal distribution, then
we have
$$W_1(\mu_X , \nu )=d_{w}\left(X\right):=\int_{-\infty}^{+\infty}\left|\mathbb{P}\left(X \leq x\right)-\Phi(x)\right| \mbox{d}x .$$
The following convergence rate in the Wasserstein-$1$ distance between $\frac{\log Z_{n}-n\mu}{\sigma\sqrt{n}}$ and the standard normal random variable is due to Grama et al. \cite{grama2017berry}. Assume that there exists a constant $\varepsilon>0 $ such that
\begin{equation}\label{a3}
	\mathbb{E}|X|^{3+\varepsilon}<\infty.
\end{equation}
Grama et al. \cite{grama2017berry} have established the following  convergence rate in the Wasserstein-$1$ distance:
\begin{equation}\label{a4}
	d_{w}\left(\frac{\log Z_{n}-n\mu}{\sigma\sqrt{n}}\right)\leq   \frac{C}{\sqrt{n}\ },
\end{equation}
where $ C $ is a positive constant.  In this paper, we are going to extend this result to the case $X$ having a moment of order $2+\delta$, $\delta \in (0, 1]$.
We get the following convergence rate in the  Wasserstein-$1$ distance for $\frac{\log Z_{n}-n\mu}{\sigma\sqrt{n}}$.  Assume that there exist  two constants $p>1$ and $\delta \in (0,1]$ such that
\begin{equation}\label{C1}
	\mathbb{E}X^{2+\delta}<\infty  \ \ \ \ \ \textrm{and} \ \ \ \ \
	\mathbb{E}\Big(\frac{Z_{1}}{m_{0}}\Big)^{p}<\infty,
\end{equation}
then it holds
\begin{equation}\label{4}
	d_{w}\left(\frac{\log Z_{n}-n\mu}{\sigma\sqrt{n}}\right)\leq   \frac{C}{n^{\delta/2}} .
\end{equation}
Moreover, the same result holds when $\frac{\log Z_{n}-n\mu}{\sigma\sqrt{n}}$ is replaced by $ \frac{n\mu-\log Z_{n}}{\sigma\sqrt{n}}$.
Notice that the last convergence rate on the right-hand side of \eqref{4} coincides with the best possible convergence rate in the Wasserstein-$1$ distance
for random walks with finite $2+\delta$ moments, and therefore inequality \eqref{4} is also the best possible.   In addition,  when $ X $ has an exponential moment and satisfies  $\displaystyle\mathbb{E} \frac{Z_{1}^{p}}{m_{0}}<\infty$, we also establish  a nonuniform Berry-Esseen bound with exponential decaying rate (see Theorem \ref{theorem2.2}).

The paper is organized as follows.  The main results are stated in Section 2.
The applications and the proofs of the main results are given
in Sections 3 and  4,  respectively.
\section{Main results}
Throughout the paper, denote by $X_i=\log m_{i},\ i\geq 0.$ Evidently, $(X_i)_{i\geq 0}$ is a sequence of i.i.d. random variables depending only on the environment $ \xi $. Let $(S_n)_{n\geq 0}$ be
the random walk associated with the branching process, which is defined as follows:
$$S_{0}=0,\ S_{n}=\log \Pi_{n}=\sum_{i=0}^{n-1}X_i,\ n\geq 1.$$ Then we have the following decomposition of $\log Z_n$, that is
\begin{equation}\label{a2}
	\log Z_{n}=S_n+\log W_{n},
\end{equation}
where $W_{n}=\frac{Z_{n}}{\Pi_{n}}. $ The normalized population size $ (W_{n})_{n\geq0} $ is a non-negative martingale under both $  \mathbb{P} $ and $ \mathbb{P}_{\xi} $, with respect to the natural filtration $(\mathcal{F}_n)_{n\geq0}$, defined by
$$
\mathcal{F}_{0}=\sigma\{\xi\}, \ \mathcal{F}_{n}=\sigma\left\{\xi, X_{k, i}, 0 \leq k \leq n-1, i \geq 1\right\}, \ n \geq 1.
$$
By Doob's martingale convergence theorem and Fatou's lemma, we can obtain that $ W_{n} $ converges a.s. to a finite limit $ W $ and $\mathbb{E}W \leq 1 $. We assume the following conditions throughout this paper:
\begin{equation}\label{a3star}
	\sigma\in(0,\infty)\quad \text{and}\quad	\mathbb{E} \frac{Z_{1}}{m_{0}} \log Z_{1}<\infty.
\end{equation}
The first condition above together with (\ref{H2}) imply in particular that 
$$Z_1\geq1\quad\text{a.s.}\quad \text{and}\quad\mathbb{P}(Z_1=1)=\mathbb{E}p_1(\xi_0)<1.$$
The second condition in (\ref{a3star}) implies that  $W_{n}$ converges to $W$ in $\mathbb{L}^{1}$ and 
$$
\mathbb{P}(W>0)=\mathbb{P}\left(Z_{n}\stackrel{n\rightarrow\infty}{\longrightarrow}\infty\right)=\lim_{n\rightarrow\infty}\mathbb{P}(Z_n>0)>0.
$$
(See for instance Tanny \cite{tanny1988necessary} and Grama et al. \cite{grama2017berry}). Therefore, it follows with the assumption (\ref{H2}) that $W>0$ and $Z_{n}\stackrel{n\rightarrow\infty}{\longrightarrow}\infty$ a.s. 

Consider the following assumptions:\\
\textbf{(A1)} There exists a constant  $\delta \in (0,1]$, such that
$$
\mathbb{E}X^{2+\delta}=\mathbb{E}\left(\log m_0\right)^{2+\delta}<\infty.
$$
\textbf{(A2)} There exists a constant $p>1$ such that
$$
\mathbb{E} \frac{m_0^{(p)}}{m_{0}^p} <\infty.$$
Under the conditions above, we obtain the following   bound for the Wasserstein-$1$ distance between $\frac{\log Z_{n}-n\mu}{\sigma\sqrt{n}}$ and the standard normal random variable.

\begin{theorem}\label{theorem2.1} Suppose that the conditions \textbf{(A1)} and \textbf{(A2)} are satisfied.  Then
	$$ d_{w}\left(\frac{\log Z_{n}-n\mu}{\sigma\sqrt{n}}\right)\leq \frac{C}{n^{\delta/2}} . $$
	Moreover, the same inequality holds when $\frac{\log Z_{n}-n\mu}{\sigma\sqrt{n}}$ is replaced by  $-\frac{\log Z_{n}-n\mu}{\sigma\sqrt{n}}.$     
\end{theorem}

Next, we consider the following stronger conditions than \textbf{(A1)} and \textbf{(A2)}:\\[8pt]
\textbf{(A3)} There exists a constant $\lambda_{0}>0$ such that
$$
\mathbb{E} e^{\lambda_{0} X}=\mathbb{E} m_{0}^{\lambda_{0}}<\infty.
$$
\textbf{(A4)} There exists a constant $p>1$ such that
$$
\mathbb{E}  \frac{Z_{1}^{p}}{m_{0}} =\mathbb{E} \frac{m_0^{(p)}}{m_0} <\infty.
$$
We have the following nonuniform Berry-Esseen bound with exponential decay rate under the conditions \textbf{(A3)} and \textbf{(A4)}, which is of independent interest. Such type of result can be found  in Fan et al. \cite{fan2017non}, where an exponential nonuniform Berry-Esseen bound for martingales has been established.
\begin{theorem}\label{theorem2.2} Suppose that the conditions  \textbf{(A3)} and \textbf{(A4)} are satisfied. Then for any  $ x \in \mathbb{R} $,
	\begin{equation}\label{2.9}
		\left|\mathbb{P}\left(\frac{\log Z_{n}-n \mu}{\sigma \sqrt{n}} \leq x\right)-\Phi(x)\right|\leq	C\frac{1}{\sqrt{n}}(1+x^{2})\exp \left\{-\frac{\hat{x}^{2}}{2}\right\},
	\end{equation}
	where $ \hat{x}=\frac{|x|}{\sqrt{1+c|x|/\sqrt{n}}}.$   
\end{theorem}

In Grama et al. \cite{grama2017berry}, Cram\'{e}r moderate deviations (cf. Theorem 1.3) have been obtained under the conditions \textbf{(A3)} and \textbf{(A4)}.
Compared to their results,  the interesting feature  of \eqref{2.9} is that it holds for any $x\in \mathbb{R}$ rather than only for  $0<|x| =o(\sqrt{n})$.

\section{Application to Interval Estimation}
In this section, we will discuss some applications of the main results to the confidence interval estimation for the criticality parameter $\mu$ and the population size $Z_n$.
\subsection{Confidence Intervals for $\mu$}
When $ \sigma $ is known,  Theorem \ref{theorem2.1} can be applied to construct the confidence interval for $\mu$.

\begin{proposition} Suppose that the conditions \textbf{(A1)} and \textbf{(A2)} are satisfied. Let $ \kappa_{n} \in (0,1) $, such that
	\begin{align}\label{k6}
		\left|\log \kappa_{n}\right|=o\left(\log n \right),\quad n\rightarrow \infty.
	\end{align}
	Then for $n$ large enough, $ \left[A_{n},\,B_{n}\right] $, with
	$$A_{n}=\frac{\log Z_{n}}{n}-\frac{\sigma\Phi^{-1}\left(1-\frac{\kappa_{n}}{2}\right)}{\sqrt{n}}  \ \ \
	\textrm{and} \ \ \
	B_{n}=\frac{\log Z_{n}}{n}+\frac{\sigma\Phi^{-1}\left(1-\frac{\kappa_{n}}{2}\right)}{\sqrt{n}}, $$
	is a $ 1-\kappa_{n} $ confidence interval for $ \mu $.  
\end{proposition}

\begin{pf} Inequality \eqref{m3.1h4} implies that
\begin{align}\label{k7}
	\frac{\mathbb{P}\left(\frac{\log Z_{n}-n\mu}{\sigma\sqrt{n}}\geq x\right)}{1-\Phi(x)}=1+o(1)\quad\text{and}\quad \frac{\mathbb{P}\left(\frac{\log Z_{n}-n\mu}{\sigma\sqrt{n}}\leq -x\right)}{\Phi(-x)}=1+o(1)
\end{align}
uniformly for $ 0\leq x=o\left(\sqrt{\log n }\right).$ For $ p\searrow 0 $, the quantile function of the standard normal distribution has the following asymptotic expansion $$ \Phi^{-1}(p)=-\sqrt{\log \frac{1}{p^{2}}-\log \log \frac{1}{p^{2}}-\log (2 \pi)}+o(1) .$$ In particular,
when $\kappa_{n}$ satisfies the condition \eqref{k6}, the upper $\left(1-\kappa_{n}/2\right)$-th quantile of standard normal distribution satisfies $  \Phi^{-1}\left(1-\frac{\kappa_{n}}{2}\right)=-\Phi^{-1}\left(\frac{\kappa_{n}}{2}\right)=O\big(\sqrt{\left|\log \kappa_{n}\right|}\big) $,  is of order $ o\left(\sqrt{\log n }\right) .$ Then, applying the last equality to \eqref{k7},  we have as $ n\rightarrow \infty $,
\begin{align}\label{h10}
	\mathbb{P}\left(\frac{\log Z_{n}-n\mu}{\sigma\sqrt{n}}\geq \Phi^{-1}\left(1-\frac{\kappa_{n}}{2}\right)\right) \sim \frac{\kappa_{n}}{2}
\end{align}
and
\begin{align}\label{h11}
	\mathbb{P}\left(\frac{\log Z_{n}-n\mu}{\sigma\sqrt{n}}\leq -\Phi^{-1}\left(1-\frac{\kappa_{n}}{2}\right)\right) \sim \frac{\kappa_{n}}{2}.
\end{align}
Therefore, as $ n\rightarrow \infty $,
\begin{align}\label{eq10*}
	\mathbb{P}\left(-\Phi^{-1}\left(1-\frac{\kappa_{n}}{2}\right)\leq\frac{\log Z_{n}-n\mu}{\sigma\sqrt{n}}\leq\Phi^{-1}\left(1-\frac{\kappa_{n}}{2}\right)\right)\sim 1-\kappa_n;
\end{align}
which implies $\mu\in[A_{n},\,B_{n}]$ with probability $1-\kappa_n$ for $n$ large enough.\hfill$\Box$
\end{pf}

\subsection{Confidence Intervals for $ Z_{n} $}
When parameters $ \mu $ and $ \sigma $ are known, we can also apply Theorem \ref{theorem2.1}  and Theorem \ref{theorem2.2} to construct confidence intervals for $ Z_{n} $. Such type of results can be used to predict the future population size $Z_n$.
\begin{proposition}\label{Proposition3.3}	Let $ \kappa_{n} \in (0,1) $. Consider the following two groups of conditions:
	\begin{itemize}
		\item[\textbf{(B1)}]
		Suppose that the conditions \textbf{(A1)} and \textbf{(A2)} are satisfied,  and that
		\begin{eqnarray}\label{kk6}
			\left|\log \kappa_{n}\right|=o\left(\log n \right), \ \ \  n\rightarrow \infty .
		\end{eqnarray}
		\item[\textbf{(B2)}]
		Suppose that the conditions \textbf{(A3)} and \textbf{(A4)} are satisfied, and that
		\begin{eqnarray}\label{m12}
			\left|\log \kappa_{n}\right|=o\big(n^{1/3} \big), \ \ \  n\rightarrow\infty.
		\end{eqnarray}
	\end{itemize}
	If $\textbf{(B1)}$ or $\textbf{(B2)}$ holds, then  for $n$ large enough, $ \left[A_{n},B_{n}\right] $ is the confidence interval of $ Z_{n} $ with confidence level $ 1-\kappa_{n} $, where
	$$A_{n}=\exp\left\{n\mu-\sigma\sqrt{n}\Phi^{-1}\left(1-\frac{\kappa_{n}}{2}\right)\right\}\quad \text{and}\quad B_{n}=\exp\left\{n\mu+\sigma\sqrt{n}\Phi^{-1}\left(1-\frac{\kappa_{n}}{2}\right)\right\}.$$	 
\end{proposition}

\begin{pf}
Assume $\textbf{(B1)}$. By inequality \eqref{m3.1h4} and \eqref{eq10*}, we have
\begin{align*}
	&\mathbb{P}\left(-\Phi^{-1}\left(1-\frac{\kappa_{n}}{2}\right)\leq\frac{\log Z_{n}-n\mu}{\sigma\sqrt{n}}\leq\Phi^{-1}\left(1-\frac{\kappa_{n}}{2}\right)\right)\notag\\
	=\,&\mathbb{P}\bigg(\exp\left\{n\mu-\sigma\sqrt{n}\Phi^{-1}\left(1-\frac{\kappa_{n}}{2}\right)\right\}\leq Z_{n}\leq\exp\left\{n\mu+\sigma\sqrt{n}\Phi^{-1}\left(1-\frac{\kappa_{n}}{2}\right)\right\}\bigg)\sim 1-\kappa_n,
\end{align*}
as $n\rightarrow \infty$. From the equality above, when $n$ is large enough, $ \left[A_{n},B_{n}\right] $ is the confidence interval for $Z_{n}$ with confidence level $ 1-\kappa_{n} $.

Now, assume $\textbf{(B2)}$. By Theorem \ref{theorem2.2}, we have
\begin{align}\label{m13}
	\frac{\mathbb{P}\left(\frac{\log Z_{n}-n\mu}{\sigma\sqrt{n}}\geq x\right)}{1-\Phi(x)}=1+o(1)\quad\text{and}\quad \frac{\mathbb{P}\left(\frac{\log Z_{n}-n\mu}{\sigma\sqrt{n}}\leq -x\right)}{\Phi(-x)}=1+o(1)
\end{align}
uniformly for $ 0\leq x=o(n^{1/6}),$ as $ n\rightarrow\infty $.
When $\kappa_{n}$ satisfies the condition \eqref{m12}, the upper $\left(1-\kappa_{n} / 2\right)$-th  quantile of the standard normal distribution satisfies $\Phi^{-1}\left(1-\kappa_{n} / 2\right)=-\Phi^{-1}\left(\kappa_{n} / 2\right)=O\left(\sqrt{\left|\log \kappa_{n}\right|}\right)$, which is of order $o\left(n^{1/6}\right)$. Then by \eqref{m13}, we get
\begin{equation*}
	\mathbb{P}\left(-\Phi^{-1}\left(1-\frac{\kappa_{n}}{2}\right)\leq\frac{\log Z_{n}-n\mu}{\sigma\sqrt{n}}\leq\Phi^{-1}\left(1-\frac{\kappa_{n}}{2}\right)\right)\sim 1-\kappa_n,
\end{equation*}
as $n\rightarrow\infty$. From the equality above, the result still holds.
\hfill$\Box$
\end{pf}

\section{Proof of Main Results}
In this section, we give the proofs of Theorems \ref{theorem2.1} and \ref{theorem2.2} respectively.
\subsection{Preliminary Lemmas for Theorem \ref{theorem2.1}}
For the sum of independent random variables, Bikelis \cite{bikelis1966estimates} gives the following nonuniform Berry-Esseen bound.
\begin{lemma}\label{lemma4.1} Let $Y_{1}, \ldots, Y_{n}$ be independent random variables satisfying $\mathbb{E} Y_{i}=0$ and  $\mathbb{E}\left|Y_{i}\right|^{2+\delta}<\infty$, $i=1, \ldots, n$, for some positive constant $\delta \in(0,1] .$ Assume that $\sum\limits_{i=1}^{n}\mathbb{E}Y_{i}^{2}=1$. Then for any $x\in\mathbb{R}$,
	$$\left|\mathbb{P}\bigg(\sum_{i=1}^{n}Y_{i}\leq x\bigg)-\Phi(x)\right|\leq \frac{C}{1+|x|^{2+\delta}}\sum\limits_{i=1}^{n} \mathbb{E}\left|Y_{i}\right|^{2+\delta}.$$    
\end{lemma}

For the convergence rate of $W_n$, we also have the following result. See Lemmas 2.3 and 2.4  in Grama et al. \cite{grama2017berry}.
See also \cite{fan2020uniform}.
\begin{lemma}\label{lemma4.3} Assume \textbf{(A1)} and \textbf{(A2)}. Then for any $\ q\in(1, 1+\delta)$,
	\begin{equation}\label{eqLiu3}
		\mathbb{E}|\log W|^q<\infty
	\end{equation}
	and
	\begin{equation}\label{eqLiu4}
		\sup_{n\in\mathbb{N}}\mathbb{E}|\log W_n|^q<\infty.
\end{equation}    
\end{lemma}
\begin{lemma} \label{lemma4.4}	Assume \textbf{(A1)} and \textbf{(A2)}. Then there exists constant $\gamma\in(0,1)$, such that for any $n\geq0$,
	$$\mathbb{E}|\log W_n-\log W|\leq C\, \gamma^n.$$		  
\end{lemma}

Using the lemmas above, we obtain the following lemma, which is a refinement of Lemma 2.5 in Grama et al. \cite{grama2017berry}.
\begin{lemma}\label{lemma4.5}  Assume \textbf{(A1)} and \textbf{(A2)}. Then   for any   $ x \in \mathbb{R} $ and $\delta^{'}\in(0,\delta)$,
	\begin{equation}\label{7}
		\mathbb{P}\left(\frac{\log Z_{n}-n \mu}{\sigma \sqrt{n}} \leq x,  \frac{S_{n}-n \mu}{\sigma \sqrt{n}} \geq x\right) \leq C\frac{1}{n^{\delta / 2}}\frac{1}{1+|x|^{1+\delta^{'}}}
	\end{equation}
	and
	\begin{equation}\label{11}
		\mathbb{P}\left(\frac{\log Z_{n}-n \mu}{\sigma \sqrt{n}} \geq x,  \frac{S_{n}-n \mu}{\sigma \sqrt{n}} \leq x\right) \leq
		C\frac{1}{n^{\delta / 2}}\frac{1}{1+|x|^{1+\delta^{'}}}.
\end{equation}		   
\end{lemma} 
\begin{pf}  We only prove the inequality \eqref{7}, as \eqref{11} can be proved in the same way.
When $x>\frac{\sqrt{n} \mu}{\sigma},$ we have
\begin{eqnarray}
	\mathbb{P}\left(\frac{\log Z_{n}-n \mu}{\sigma \sqrt{n}} \leq x,  \frac{S_{n}-n \mu}{\sigma \sqrt{n}} \geq x\right)  \leq  \mathbb{P}\left(  \frac{S_{n}-n \mu}{\sigma \sqrt{n}} \geq x\right) , \nonumber
\end{eqnarray}
which implies \eqref{7} by Lemma \ref{lemma4.1}. When $x < -\frac{\sqrt{n} \mu}{\sigma},$  it follows from the fact that $Z_{n}\geq 1 \ \text{a.s.}$,
$$ \frac{\log Z_{n}-n \mu}{\sigma \sqrt{n}} \geq  -\frac{\sqrt{n} \mu}{\sigma}.$$
Then for any $ x < -\frac{\sqrt{n} \mu}{\sigma} $,
\begin{eqnarray}
	\mathbb{P}\left(\frac{\log Z_{n}-n \mu}{\sigma \sqrt{n}} \leq x,  \frac{S_{n}-n \mu}{\sigma \sqrt{n}} \geq x\right) =0.\nonumber
\end{eqnarray}
Thus for any  $|x| >\frac{\sqrt{n} \mu}{\sigma}$, \eqref{7} holds.

Next, we prove that for any $|x| \leq\frac{\sqrt{n} \mu}{\sigma}$, \eqref{7} holds.
Let's consider the following notations
$$Y_{m, n}=\sum_{i=m+1}^{n} \frac{X_{i}-\mu}{\sigma \sqrt{n}},\quad  Y_{n}=Y_{0, n},\quad V_{m}=\frac{\log W_{m}}{\sigma \sqrt{n}}\ \text { and } \ D_{m, n}=V_{n}-V_{m}.$$
Set $ \alpha_{n}=\frac{1}{\sqrt{n}}$ and $m=[\sqrt{n}]$, where $ [t]$ stands for the largest integer less than $t$. From \eqref{a2},
we get
\begin{align}\label{10}
	\mathbb{P}\left(\frac{\log Z_{n}-n \mu}{\sigma \sqrt{n}} \leq x,  \frac{S_{n}-n \mu}{\sigma \sqrt{n}} \geq x\right)
	\leq\mathbb{P}\left(Y_{n}+V_{m} \leq x+\alpha_{n},  Y_{n} \geq x\right)\notag\nonumber\\
	+\ \mathbb{P}\left(\left|D_{m,  n}\right|>\alpha_{n}\right).
\end{align}
For the second term on the right-hand side of \eqref{10}, using Markov's inequality and
Lemma \ref{lemma4.4},
for any $|x|\leq\frac{\sqrt{n} \mu}{\sigma}$, there exists a constant $ \gamma \in (0,1)$ such that
\begin{align}\label{a3.4}
	\mathbb{P}\left(\left|D_{m,  n}\right|>\alpha_{n}\right) \leq \displaystyle \frac{\mathbb{E}\left|D_{m,  n}\right|}{\alpha_{n}}=& \frac{\sqrt{n}\, \mathbb{E}\left|\log W_{n}-\log W_{m}\right|}{\sigma \sqrt{n}}\leq \displaystyle  C \gamma ^{m} \nonumber  \\
	\leq&\, C \frac{1}{n^{\delta / 2}}\frac{1}{1+|x|^{2+\delta}}. \ \ \ \ \ \ \ \ \ \ \ \ \ \ \ \ \ \ \ \ \ \ \
\end{align}
The last line is obtained from the fact that $n^{\delta/2}(1+ n^{(1+\delta)/2})\gamma^m=o(1)$,   $n\rightarrow\infty$.\
For the first term on the right-hand side of \eqref{10}, we use a decomposition to dominate it. Let $$G_{m, n}(x)=\mathbb{P}\left(Y_{m, n} \leq x\right),\ G_{ n}(x)=\mathbb{P}\left(Y_{n} \leq x\right) \ \text{ and } \ v_{m}(\mbox{d} s, \mbox{d} t)=\mathbb{P}\left(Y_{m} \in \mbox{d} s, V_{m} \in \mbox{d} t\right).$$
Since $ Y_{m, n} $ and $(Y_{m},V_{m})  $ are independent, we have
\begin{eqnarray}\label{a3.5}
	&&\mathbb{P}\left(Y_{n}+V_{m} \leq x+\alpha_{n},  Y_{n} \geq x\right)\notag \\
	&&=\mathbb{P}\left(Y_{m,  n}+Y_{m}+V_{m} \leq x+\alpha_{n},  Y_{m,  n}+Y_{m} \geq x\right) \notag\\
	&&= \int\!\!\int \mathbb{P}\left(Y_{m,  n}+s+t \leq x+\alpha_{n},  Y_{m,  n}+s \geq x\right) v_{m}(\mbox{d}s,  \mbox{d}t) \notag\\
	&&=  \int\!\!\int \mathbf{1}_{\{t \leq \alpha_{n}\}}\left(G_{m,  n}\left(x-s-t+\alpha_{n}\right)-G_{m,  n}(x-s)\right) v_{m}(\mbox{d}s,  \mbox{d}t).
\end{eqnarray}
Notice that $(X_i)_{i\geq 1}$ are i.i.d., so $\displaystyle \sum_{i=m+1}^nX_i$ and $\displaystyle \sum_{i=1}^{n-m}X_i$ follow the same distribution. We have
\begin{align}\label{m3.6}
	G_{m, n}(x)
	=\mathbb{P}\left(\sum_{i=m+1}^{n} \frac{X_{i}-\mu}{\sigma \sqrt{n}} \leq x\right) =G_{n-m}\left(\frac{x \sqrt{n}}{\sqrt{n-m}}\right)=G_{n-m}\Big(x\left(1+R_{n}\right)\Big),
\end{align}
where
$R_{n}=\sqrt{\frac{n}{n-m}}\,-1$
satisfying $0 \leq R_{n} \leq \frac{C}{\sqrt{n}}.$
By the mean value theorem, we deduce  that for any  $ x\in \mathbb{R} $,
\begin{eqnarray}\label{4.45}
	\left|\Phi\left(x\left(1+R_{n}\right)\right)-\Phi(x)\right| &\leq& R_{n}\left|x\frac{1}{\sqrt{2\pi}}\exp\left\{-\frac{x^2}{2}\right\}\right| \ \leq \ C\frac{1}{\sqrt{n}}\left|x\exp\left\{-\frac{x^{2}}{2}\right\}\right|\notag\\
	&\leq& C\frac{1}{n^{\delta / 2}}\frac{1}{1+|x|^{2+\delta}}.
\end{eqnarray}
It can be obtained from Lemma \ref{lemma4.1} that
\begin{align}\label{4.46}
	\left| G_{n-m}\left(x\left(1+R_{n}\right)\right)-\Phi\left(x\left(1+R_{n}\right)\right) \right|
	&=\left|\mathbb{P}\Big(\sum_{i=1}^{n-m} \frac{X_{i}-\mu}{\sigma \sqrt{n-m}} \leq \frac{x\sqrt{n}}{\sqrt{n-m}}\Big)-\Phi\Big(\frac{x\sqrt{n}}{\sqrt{n-m}}\Big)\right|\notag\\
	&\leq C\frac{1}{(n-m)^{\delta/2}}\frac{1}{1+|x\left(1+R_{n}\right)|^{2+\delta}}\notag\\
	&\leq  C\frac{1}{n^{\delta / 2}}\frac{1}{1+|x|^{2+\delta}}.
\end{align}
Combining \eqref{m3.6}-\eqref{4.46} together, we get
$$
|G_{m,n}(x)-\Phi(x)|\leq C\frac{1}{n^{\delta / 2}}\frac{1}{1+|x|^{2+\delta}}.
$$
Therefore, by \eqref{a3.5}, we have
\begin{eqnarray}\label{u4.34}
	\mathbb{P}\left(Y_{n}+V_{m} \leq x+\alpha_{n},  Y_{n} \geq x\right)\leq J_{1}+J_{2}+J_{3},
\end{eqnarray}
where
\begin{eqnarray*}
	&&J_{1}=\int\!\!\int \mathbf{1}_{\{t \leq \alpha_{n}\}}\left|\Phi\left(x-s-t+\alpha_{n}\right)-\Phi(x-s)\right| v_{m}(\mbox{d} s,  \mbox{d} t),\\
	&&J_{2}=C\frac{1}{n^{\delta / 2}}\int\!\!\int\mathbf{1}_{\{t \leq \alpha_{n}\}}\frac{1}{1+|x-s|^{2+\delta}} v_{m}(\mbox{d} s,\mbox{d}t)
\end{eqnarray*}
and
\begin{eqnarray*}
	J_{3}=C\frac{1}{n^{\delta / 2}}\int\!\!\int\mathbf{1}_{\{t \leq \alpha_{n}\}}\frac{1}{1+|x-s-t+\alpha_{n}|^{2+\delta}}\, v_{m}(\mbox{d} s,\mbox{d}t).
\end{eqnarray*}
For $J_{1}$, by the mean value theorem, it holds
\begin{align*}
	&\mathbf{1}_{\{t \leq \alpha_{n}\}}\,\left|\Phi\left(x-s-t+\alpha_{n}\right)-\Phi(x-s)\right| \\
	\leq&\, C\,|-t+\alpha_{n}|\exp\left\{-\frac{x^{2}}{8}\right\} + \,  |-t+\alpha_{n}|\, \mathbf{1}_{\{|s| \geq 1+  \frac{1}{4 }|x| \}} +   \,  |-t+\alpha_{n}|\, \mathbf{1}_{\{|t| \geq 1+  \frac{1}{4 }|x| \}},
\end{align*}
and thus we have
\begin{equation}\label{4.37}
	J_{1}\leq J_{11}+J_{12}+J_{13},
\end{equation}
where
\begin{eqnarray}
	&&J_{11}= C\int\!\!\int |-t+\alpha_{n}|\exp\left\{-\frac{x^{2}}{8}\right\} v_{m}(\mbox{d} s,  \mbox{d} t), \nonumber \\
	&&J_{12}=  \int\!\!\int |-t+\alpha_{n}|\mathbf{1}_{\{|s| \geq 1+  \frac{1}{4 }|x| \}}  v_{m}(\mbox{d} s,  \mbox{d} t)\nonumber
\end{eqnarray}
and
\begin{eqnarray}
	J_{13}= \int\!\!\int |-t+\alpha_{n}|\mathbf{1}_{\{|t| \geq 1+  \frac{1}{4 }|x| \}} v_{m}(\mbox{d} s,  \mbox{d} t). \nonumber
\end{eqnarray}
Firstly, considering $J_{11}$, it follows from Lemma \ref{lemma4.3} that for any $x \in\mathbb{R}$,
\begin{eqnarray}
	J_{11}& \leq &  C \exp\left\{-\frac{x^{2}}{8}\right\}  \bigg(  \mathbb{E}| V_{m}| +  \alpha_{n} \bigg)\nonumber \\
	&\leq&    C\frac{1}{n^{\delta/ 2}} \frac{1}{1+|x|^{2+\delta }}.
\end{eqnarray}
For $J_{12}$, Let $\tau,\,\iota>1$, and satisfy $\frac{1}{\tau}+\frac{1}{\iota}=1$.
By H\"{o}lder's inequality, we have the following estimate: for any $x \in\mathbb{R}$,
\begin{eqnarray}\label{4.39}
	J_{12} & \leq& \iint|t| \mathbf1_{\left\{|s| \geq 1+\frac{1}{4} |x|\right\}} v_{m}(\mbox{d} s, \mbox{d} t)+\alpha_{n} \int \mathbf1_{\left\{|s| \geq 1+\frac{1}{4}|x|\right\}} v_{m}(\mbox{d} s) \nonumber\\
	&\leq&\left(\int|t|^{\tau} v_{m} (\mbox{d} t)\right)^{\frac{1}{\tau}} \left(\int \mathbf{1}_{\{|s|\geq 1+\frac{1}{4}|x|\}}  v_{m}(\mbox{d} s)\right)^{\frac{1}{\iota}}+\frac{1}{\sqrt{n}} \; \mathbb{P}\left(\left|Y_{m}\right| \geq 1+\frac{1}{4}|x|\right) \nonumber\\
	& \leq &C\frac{1}{\sqrt{n}}\left[\mathbb{P}\left(\left|Y_{m}\right| \geq 1+\frac{1}{4}|x|\right)\right]^{\frac{1}{l}}+\frac{1}{\sqrt{n}} \mathbb{P}\left(\left|Y_{m}\right| \geq 1+\frac{1}{4}|x|\right).
\end{eqnarray}
From Lemma \ref{lemma4.1}, we get
\begin{eqnarray}\label{h4.40}
	\mathbb{P}\bigg( |Y_{m} | \geq 1+  \frac{1}{4 }|x| \bigg)
	&=& \mathbb{P}\bigg( Y_{m}  \geq 1+  \frac{1}{4 }|x| \bigg)  +\mathbb{P}\bigg( Y_{m}  \leq -1-  \frac{1}{4 }|x| \bigg)  \nonumber \\
	&\leq& 1-\Phi\left(\big(1+  \frac{1}{4 }|x|\big)\frac{\sqrt{n}}{\sqrt{m}}\right)+\Phi\left(-\big(1+\frac{1}{4 }|x|\big)\frac{\sqrt{n}}{\sqrt{m}}\right)\nonumber \\
	&& +\, C\, \frac{1}{1+\left((1+|x|/4)\sqrt{n}/\sqrt{m}\right)^{2+\delta}}\sum_{i=1}^{m}\mathbb{E}\left|\frac{X_{i}-\mu}{\sigma\sqrt{m}}\right|^{2+\delta}\nonumber \\
	&\leq&  C\frac{1}{n^{(1+\delta)/ 2}}\frac{1}{1+|x|^{2+\delta}}.
\end{eqnarray}
If $\delta^{'}=\frac{3+2\delta}{\iota}-\delta$,  we can obtain
\begin{equation}\label{defD}
	\delta^{'}=3+\delta-\frac{3+2\delta}{\tau}.
\end{equation}
Substituting \eqref{h4.40} into \eqref{4.39},  we get,  for any $|x|\leq\mu\sqrt{n}/\sigma$,
\begin{eqnarray}
	J_{12}\leq C\frac{1}{n^{\delta / 2}}\frac{1}{1+|x|^{1+\delta^{'}}}.
\end{eqnarray}
Next, we turn to $J_{13}$. Set $\displaystyle p=\frac{\delta^{'}+\delta}{2}$. Using Markov's inequality and Lemma \ref{lemma4.3}, we have for any $|x| \leq \frac{\sqrt{n} \mu}{\sigma}$,
\begin{eqnarray}\label{4.40}
	J_{13}
	& \leq &
	\int |t|\mathbf{1}_{\{|t| \geq 1+  \frac{1}{4 }|x| \}} v_{m}(  \mbox{d} t)+\int \alpha_{n}\mathbf{1}_{\{|t| \geq 1+  \frac{1}{4 }|x| \}} v_{m}( \mbox{d} t)\nonumber \\
	&\leq& \int \frac{|t|^{p}}{(1+|x|/ 4)^{p}}\;|t|\; v_{m}(\mbox{d} t)+\frac{1}{\sqrt{n}} \mathbb{P}\left(\left|V_{m}\right| \geq 1+\frac{1}{4}|x|\right)\nonumber\\
	&\leq& \frac{1}{(1+|x| / 4)^{p}}\, \mathbb{E}\left|V_{m}\right|^{1+p}+\frac{1}{\sqrt{n}} \frac{1}{(1+|x| / 4)^{1+p}}\,\mathbb{E}\left|V_{m}\right|^{1+p}\nonumber\\
	&\leq& C\frac{1}{n^{\delta / 2}}\frac{1}{1+|x|^{1+\delta^{'}}}.
\end{eqnarray}
Combining \eqref{4.37}-\eqref{4.40}, we get
\begin{eqnarray}\label{4.41}
	J_{1}\leq   C\frac{1}{n^{\delta/ 2}} \frac{1}{1+|x|^{1+\delta^{'}}}.
\end{eqnarray}
	Now we discuss the upper bounds for $J_{2}$ and $J_{3}$ respectively. By an argument similar to that of \eqref{h4.40}, we can obtain for any $x\in\mathbb{R}$,
	\begin{eqnarray}\label{u4.39}
		J_{2}
		&=&C\frac{1}{n^{\delta / 2}}\int\!\!\int\mathbf{1}_{\{t \leq \alpha_{n}\}}\frac{1}{1+|x-s|^{2+\delta}} v_{m}(\mbox{d} s, \mbox{d}t)\notag \\
		&\leq&   C\frac{1}{n^{\delta / 2}}\left(\int_{|s|< 1+|x|/2}\frac{1}{1+|x-s|^{2+\delta}} v_{m}(\mbox{d} s)+\int_{|s| \geq 1+ |x|/2}\frac{1}{1+|x-s|^{2+\delta}} v_{m}(\mbox{d} s)\right)  \notag \\
		&\leq&    C\frac{1}{n^{\delta / 2}}\left[  \frac{1}{1+|x/2|^{2+\delta}}+\mathbb{P}\bigg( |Y_{m} | \geq 1+  \frac{1}{2}|x| \bigg)\right]\notag  \\
		&\leq&   C\frac{1}{n^{\delta / 2}}\frac{1}{1+|x|^{2+\delta}}.
	\end{eqnarray}
	For $J_{3}$, using an argument similar to that of \eqref{h4.40} and \eqref{4.40},  we have  for any $|x| \leq   \frac{\sqrt{n} \mu}{\sigma},$
	\begin{eqnarray}\label{fgsg45}
		J_{3}
		&=&  C\frac{1}{n^{\delta / 2}}\int\!\!\int \mathbf{1}_{\{t \leq \alpha_{n}\}} \frac{1}{1+|x-s-t|^{2+\delta}} v_{m}(\mbox{d} s,  \mbox{d} t)\nonumber  \\
		&\leq& C\frac{1}{n^{\delta / 2}}\left(\int\!\!\int_{|s+t|< 2+ |x|/2}\frac{1}{1+|x/2|^{2+\delta}} v_{m}(\mbox{d} s,  \mbox{d} t)\right. \nonumber \\
		&&\left. +\int_{|s| \geq 1+ |x|/4} v_{m}(\mbox{d} s)+  \int_{|t| \geq 1+ |x|/4} v_{m}( \mbox{d} t) \right)  \notag \\
		&\leq&    C\frac{1}{n^{\delta / 2}}\left[  \frac{1}{1+|x/2|^{2+\delta}}+\mathbb{P}\left( |Y_{m} |  \geq 1+  \frac{1}{4}|x|\right)+  \mathbb{P}\bigg( | V_{m}| \geq 1+  \frac{1}{4}|x| \bigg)\right] \notag  \\
		&\leq&   C\frac{1}{n^{\delta / 2}}\frac{1}{1+|x|^{2+\delta}}.
	\end{eqnarray}
	Then, substituting \eqref{4.41}-\eqref{fgsg45} into \eqref{u4.34}, we get
	\begin{equation}\label{4.44}
		\mathbb{P}\left(Y_{n}+V_{m} \leq x+\alpha_{n},  Y_{n} \geq x\right)\leq   C\frac{1}{n^{\delta / 2}}\frac{1}{1+|x|^{1+\delta^{'}}}.
	\end{equation}
	Finally, from \eqref{10}, \eqref{a3.4} and \eqref{4.44}, we have \eqref{7} holds for any $|x|\leq \frac{\sqrt{n} \mu}{\sigma}.$
	\hfill$\Box$
	\end{pf}
	\subsection{Proof of Theorem \ref{theorem2.1}}
	From \eqref{a2}, it holds
	$$
	\frac{\log Z_{n}-n \mu}{\sigma \sqrt{n}}=\frac{S_{n}-n \mu}{\sigma \sqrt{n}}+\frac{\log W_{n}}{\sigma \sqrt{n}}.
	$$
	By Lemma \ref{lemma4.1} and \textbf{(A1)}, we have
	\begin{eqnarray}\label{44}
		\left|\mathbb{P}\left(\frac{S_{n}-n \mu}{\sigma \sqrt{n}} \leq x\right)-\Phi(x)\right|
		&\leq& C \sum_{i=1}^{n} \mathbb{E}\left|\frac{X_{i}-\mu}{\sigma \sqrt{n}}\right|^{2+\delta} \frac{1}{(1+|x|)^{2+\delta}} \nonumber  \\
		&\leq& C \frac{1}{n^{\delta / 2}} \frac{1}{1+|x|^{2+\delta}}.
	\end{eqnarray}
	Notice that
	\begin{eqnarray}
		&&\mathbb{P}\left(\frac{\log Z_{n}-n \mu}{\sigma \sqrt{n}} \leq x\right) \notag\\
		&=&\,  \mathbb{P}\left(\frac{\log Z_{n}-n \mu}{\sigma \sqrt{n}} \leq x,  \frac{S_{n}-n \mu}{\sigma \sqrt{n}} \leq x\right)+\mathbb{P}\left(\frac{\log Z_{n}-n \mu}{\sigma \sqrt{n}} \leq x,  \frac{S_{n}-n \mu}{\sigma \sqrt{n}}>x\right) \notag\\
		&=&\, \mathbb{P}\left(\frac{S_{n}-n \mu}{\sigma \sqrt{n}} \leq x\right)+\mathbb{P}\left(\frac{\log Z_{n}-n \mu}{\sigma \sqrt{n}} \leq x,  \frac{S_{n}-n \mu}{\sigma \sqrt{n}}>x\right) \notag \\
		&& \qquad -\ \mathbb{P}\left(\frac{\log Z_{n}-n \mu}{\sigma \sqrt{n}}>x,  \frac{S_{n}-n \mu}{\sigma \sqrt{n}} \leq x\right).
	\end{eqnarray}
	Applying Lemma \ref{lemma4.5} to the last equality, we get for any $x \in \mathbb{R},$
	\begin{eqnarray}\label{m3.14}
		\left|\mathbb{P}\left(\frac{\log Z_{n}-n \mu}{\sigma \sqrt{n}} \leq x\right)-\mathbb{P}\left(\frac{S_{n}-n \mu}{\sigma \sqrt{n}} \leq x\right)\right|\leq C\frac{1}{n^{\delta / 2}}\frac{1}{1+|x|^{1+\delta^{'}}}.
	\end{eqnarray}
	Combining \eqref{44} and \eqref{m3.14}, we have for any $x \in \mathbb{R},$
	\begin{eqnarray}
		\left|\mathbb{P}\left(\frac{\log Z_{n}-n \mu}{\sigma \sqrt{n}} \leq x\right)-\Phi(x)\right|
		&\leq&\left|\mathbb{P}\left(\frac{\log Z_{n}-n \mu}{\sigma \sqrt{n}} \leq x\right)-\mathbb{P}\left(\frac{S_{n}-n \mu}{\sigma \sqrt{n}} \leq x\right)\right| \nonumber\\
		&&+\left|\mathbb{P}\left(\frac{S_{n}-n \mu}{\sigma \sqrt{n}} \leq x\right)-\Phi(x)\right| \nonumber\\
		& \leq& C\frac{1}{n^{\delta / 2}}\frac{1}{1+|x|^{1+\delta^{'}}}.\label{m3.1h4}
	\end{eqnarray}
	Thus
	\begin{align*}
		d_{w}\left(\frac{\log Z_{n}-n\mu}{\sigma\sqrt{n}}\right)
		&= \int_{-\infty}^{\infty} \left|\mathbb{P}\left(\frac{\log Z_{n}-n \mu}{\sigma \sqrt{n}} \leq x\right)-\Phi(x)\right|\mbox{d}x\\
		&\leq \frac{C}{n^{\delta / 2}}\int_{-\infty}^{\infty} \frac{1}{1+|x|^{1+\delta^{'}}}\mbox{d}x\nonumber \\
		&\leq   \frac{C}{n^{\delta/2}} ,
	\end{align*}
	which gives the first desired inequality of Theorem \ref{theorem2.1}. By a similar argument,
	it is easy to show that the same result holds when $\frac{\log Z_{n}-n\mu}{\sigma\sqrt{n}}$ is replaced by  $-\frac{\log Z_{n}-n\mu}{\sigma\sqrt{n}}.$
	\hfill$\Box$
\subsection{Preliminary Lemmas for Theorem \ref{theorem2.2}}
To prove Theorem \ref{theorem2.2}, we shall make use of the following lemma (see Theorem 3.1 in Grama et al. \cite{grama2017berry}). The lemma shows that conditions (\textbf{A3}) and \textbf{(A4)} imply the existence of a harmonic moment of positive order $ a_{0} $.
\begin{lemma}\label{lemma4.6}  Assume \textbf{(A3)} and \textbf{(A4)}. Then there exists a positive $a_{0}$ such that
	\begin{equation}\label{eqLidsu3}
		\mathbb{E} W ^{-a_{0}}<\infty \ \ \ \ \  and \ \ \ \ \sup_{n\in\mathbb{N}}\mathbb{E}  W_n^{-a_{0}}<\infty.
\end{equation} 
\end{lemma}

The next lemma shows that in the case of i.i.d., Cram\'{e}r's condition (\textbf{A3}) and Bernstein's condition \textbf{(A3$'$)}  are equivalent. See Fan et al. \cite{FGL13}.
\begin{lemma} \label{lemma4.7}  Condition \textbf{(A3)} is  equivalent to the following condition
	\textbf{(A3$'$)}. There exists a constant $H>0$  such that for any  $k \geq 2$,
	\begin{equation}\label{eq38}
		\mathbb{E}(X-\mu)^{k} \leq \frac{1}{2} k ! H^{k-2} \mathbb{E}(X-\mu)^{2}.
\end{equation}  
\end{lemma}

The following lemma gives two Bernstein type inequalities for $ \log Z_{n} $.
\begin{lemma}\label{lemma4.8}	Suppose that the conditions \textbf{(A3)} and \textbf{(A4)} are satisfied. Then for any  $x\geq0$,
	\begin{equation}\label{3.17}
		\mathbb{P}\left(\frac{\log Z_{n}-n \mu}{\sigma \sqrt{n}} \geq x\right) \leq 2 \exp \left\{-\frac{x^{2}}{2\left(1+cx/\sqrt{n}\right)}\right\}
	\end{equation}	
	and
	\begin{equation}\label{3.18}
		\mathbb{P}\left(\frac{\log Z_{n}-n \mu}{\sigma \sqrt{n}} \leq -x\right) \leq C \exp \left\{-\frac{x^{2}}{2\left(1+cx/\sqrt{n}\right)}\right\},
	\end{equation}
	where $c, C$ are two positive  constants.
\end{lemma}
\begin{pf}
Since the Cram\'{e}r's condition \textbf{(A3)} is equivalent to the Bernstein condition \textbf{(A3$'$)}, we only need to prove Lemma \ref{lemma4.8} under the conditions \textbf{(A3$'$)} and \textbf{(A4)}.

We first give a proof for \eqref{3.17}.
Denote
$$
\eta_{n, i}=\frac{X_{i}-\mu}{\sigma \sqrt{n}}, \ i=1, \dots, n.
$$
It's easy to see that
$$
\frac{\log Z_{n}-n \mu}{\sigma \sqrt{n}}=\sum_{i=1}^{n}\eta_{n,i}+\frac{\log W_{n}}{\sigma \sqrt{n}},
$$
where $ \sum\limits_{i=1}^{n} \eta_{n,i} $ is a sum of i.i.d. random variables.
Then we have for any $ x\geq 0 $,
\begin{eqnarray}\label{c3.16}
	\mathbb{P}\left(\frac{\log Z_{n}-n \mu}{\sigma \sqrt{n}} \geq x\right)&=& \mathbb{P}\left(\sum_{i=1}^{n}\eta_{n, i}+\frac{\log W_{n}}{\sigma \sqrt{n}} \geq x\right)\notag\\
	&\leq&\,   I_{1}+I_{2},
\end{eqnarray}
where $$I_{1}=\mathbb{P}\left(\sum_{i=1}^{n}\eta_{n,i} \geq x-\frac{x^{2}}{\sigma\sqrt{n}}\right)\ \ \ \text{and}\ \ \  I_{2}=\mathbb{P}\left(\frac{\log W_{n}}{\sigma \sqrt{n}} \geq \frac{x^{2}}{\sigma\sqrt{n}}\right).$$
Applying Bernstein's inequality to $I_1$, we obtain for any $x\in  [0,\,\frac{\sigma\sqrt{n}}{2} ) $,
\begin{equation}\label{c3.17}
	I_{1} \leq \exp \Bigg\{-\frac{x^{2}\left(1-\frac{x}{\sigma \sqrt{n}}\right)^{2}}{2\left(1+\frac{H}{\sigma \sqrt{n}} x\big(1-\frac{x}{\sigma \sqrt{n}}\big)\right)}\Bigg\}\leq \exp
	\Bigg\{-\frac{x^{2}}{2\left(1+cx/\sqrt{n}\right)}\Bigg\}.
\end{equation}
By Markov's inequality and the fact $ \mathbb{E} W_{n}=1 $, we have for any $x\in [0,\, \frac{\sigma\sqrt{n}}{2} )$,	
\begin{equation}\label{c3.18}
	I_{2}=\mathbb{P}\left(W_{n} \geq \exp \left\{x^{2}\right\}\right)\leq \exp \left\{-x^{2}\right\} \mathbb{E} W_{ n}=\exp \left\{-x^{2}\right\} .
\end{equation}
Combining \eqref{c3.16}-\eqref{c3.18} together,  we find that \eqref{3.17} holds for any  $x\in [0,\frac{\sigma\sqrt{n}}{2} )$.
When $ x\geq \frac{\sigma\sqrt{n}}{2} $, it holds that
\begin{align}\label{c3.20}
	\mathbb{P}\left(\frac{\log Z_{n}-n \mu}{\sigma \sqrt{n}} \geq x\right) \leq I_{3}+I_{4},
\end{align}
where
\begin{align*}
	I_{3}=\mathbb{P}\left(\sum_{i=1}^{n} \eta_{n,i} \geq \frac{x}{2}\right)\ \ \ \text{and} \ \ \ I_{4}=\mathbb{P}\left(\frac{\log W_{n}}{\sigma \sqrt{n}} \geq \frac{x}{2}\right).
\end{align*}
By the same arguments as  the proofs of $I_{1}$ and $I_{2}$, we have for any $ x\geq\frac{\sigma\sqrt{n}}{2} $,
\begin{align}\label{c3.21}
	I_{3}
	\leq \exp \left\{-\frac{(x / 2)^{2}}{2\left(1+\frac{H}{\sigma \sqrt{n}} \frac{x}{2}\right)}\right\}
	\leq \exp \left\{-\frac{x^{2}}{2\left(1+cx/\sqrt{n}\right)}\right\}
\end{align}
and
\begin{align}\label{c3.22}
	I_{4} 	\leq \exp \left\{-\frac{x \sigma \sqrt{n}}{2}\right\} \mathbb{E} W_{n}
	\leq \exp\left\{-\frac{x^{2}}{2\left(1+cx/\sqrt{n}\right)}\right\},
\end{align}
with $c$ large enough.
Combining \eqref{c3.20}-\eqref{c3.22} together, we get \eqref{3.17}  for  $x\geq\frac{\sigma\sqrt{n}}{2}$.

Next, we will prove \eqref{3.18}. Let $ a_{0} $ be a positive constant given by Lemma \ref{lemma4.6}. Then it holds for any $ x \geq 0,$
\begin{eqnarray}
	\mathbb{P}\left(\frac{\log Z_{n}-n \mu}{\sigma \sqrt{n}} \leq -x\right)&=& \mathbb{P}\left(-\sum_{i=1}^{n}\eta_{n,i}-\frac{\log W_{n}}{\sigma \sqrt{n}} \geq x\right) \nonumber \\
	&\leq& I_{5}+I_{6},\label{3.19}
\end{eqnarray}
where $$I_{5}=\mathbb{P}\left(-\sum_{i=1}^{n}\eta_{n,i} \geq x-\frac{x^{2}}{a_{0}\sigma\sqrt{n}}\right) \ \ \ \  \text{and}\  \ \ I_{6}=\mathbb{P}\left(-\frac{\log W_{n}}{\sigma \sqrt{n}} \geq \frac{x^{2}}{a_{0}\sigma\sqrt{n}}\right).$$
The upper bounds for $ I_{5} $ and $ I_{6} $ are given respectively as follows. Using Bernstein's inequality, we get for any $ x\in [0,\,\frac{a_{0}\sigma\sqrt{n}}{2} ) $,
\begin{equation}\label{3.20}
	I_{5}\leq \exp \Bigg\{-\frac{x^{2}\left(1-\frac{x}{a_{0}\sigma \sqrt{n}}\right)^{2}}{2\Big(1+\frac{H}{\sigma \sqrt{n}} x\big(1-\frac{x}{a_{0}\sigma \sqrt{n}}\big)\Big)}\Bigg\}\leq \exp \left\{-\frac{x^{2}}{2\left(1+cx/\sqrt{n}\right)}\right\}.
\end{equation}
And by Markov's inequality, we get for any $x\in  [0,\,\frac{a_{0}\sigma\sqrt{n}}{2} ),$
\begin{eqnarray}\label{i4.62}
	I_{6} &=&\mathbb{P}\left(\log W_{n} \leq -\frac{x^{2}}{a_{0}}\right)=\mathbb{P}\Big(W_n^{-a_{0}}\geq\exp{\left\{x^2\right\}}\Big)\nonumber \\
	&\leq& \exp\{-x^{2}\}\mathbb{E}W_{n}^{-a_{0}}<\infty.
\end{eqnarray}
Combining \eqref{3.19}-\eqref{i4.62} together, we get \eqref{3.18}  for any $x\in  [0,\, \frac{a_{0}\sigma\sqrt{n}}{2} )$.
When $ x\geq\frac{a_{0}\sigma\sqrt{n}}{2} $,  it holds that
\begin{equation*}
	\mathbb{P}\left(\frac{\log Z_{n}-n \mu}{\sigma \sqrt{n}} \leq -x\right) \leq I_{7}+I_{8},
\end{equation*}
where $$I_{7}=\mathbb{P}\bigg(-\sum_{i=1}^{n} \eta_{n,i} \geq \frac{x}{2}\bigg)
\ \ \  \text{and}\ \ \  \ I_{8}=\mathbb{P}\left( -\frac{\log W_{n}}{\sigma \sqrt{n}} \geq \frac{x}{2}\right).$$
Again by Bernstein's inequality, we obtain for any  $ x\geq\frac{a_{0}\sigma\sqrt{n}}{2} $,
\begin{align*}
	I_{7} \leq \exp \bigg\{-\frac{(x / 2)^{2}}{2\big(1+\frac{H}{\sigma \sqrt{n}} \frac{x}{2}\big)}\bigg\} \leq \exp \left\{-\frac{x^{2}}{2\left(1+cx/\sqrt{n}\right)}\right\}.
\end{align*}
And by the same arguments as  the proof of $ I_{6} $, we can get for any  $ x\geq\frac{a_{0}\sigma\sqrt{n}}{2} $,
\begin{equation*}
	I_{8}=\mathbb{P}\left(\log W_{n} \leq -\frac{\sigma \sqrt{n}x}{2}\right)\leq\exp\left\{-\frac{a_{0}\sigma\sqrt{n}x}{2}\right\}\mathbb{E}W_{n}^{-a_{0}}\leq C\exp \left\{-\frac{x^{2}}{2\left(1+cx/\sqrt{n}\right)}\right\}.
\end{equation*}
This completes the proof of lemma. \hfill$\Box$
\end{pf}

The following lemma  is a direct consequence of Theorem 1.3 in \cite{grama2017berry} and its detailed proof  can be found in \cite{grama2017berry}.
\begin{lemma}\label{lemma4.9} Assume \textbf{(A3)} and \textbf{(A4)}. Then for any $x\in\left[0,\, n^{1/4}\right]$,
	\begin{equation*}
		\left|\log \frac{\mathbb{P}\left( \frac{\log Z_{n}-n \mu}{\sigma \sqrt{n}}>x \right) }{1-\Phi(x)} \right|\leq C\frac{1+x^{3}}{\sqrt{n}}
	\end{equation*}
	and
	\begin{equation*}
		\left|\log \frac{\mathbb{P}\left( \frac{\log Z_{n}-n \mu}{\sigma \sqrt{n}}\leq -x \right) }{\Phi(-x)} \right|\leq C\frac{1+x^{3}}{\sqrt{n}}.
\end{equation*}  
\end{lemma}  
	
\subsection{Proof of Theorem \ref{theorem2.2}}
When $x\in\left[ 0,\,n^{1/4}\right] $, by the first inequality in Lemma \ref{lemma4.9}, we get
\begin{align}\label{3.29}
	\mathbb{P}\left(\frac{\log Z_{n}-n \mu}{\sigma \sqrt{n}}\leq x\right)-\Phi(x)
	&=-\left[\mathbb{P}\left(\frac{\log Z_{n}-n \mu}{\sigma \sqrt{n}} > x\right)-(1-\Phi(x))\right]\notag\\
	&\geq-\left[(1-\Phi(x))\exp\left\{C\frac{1+x^{3}}{\sqrt{n}}\right\}-(1-\Phi(x))\right]\notag\\
	&=-(1-\Phi(x))\left(\exp\left\{C\frac{1+x^{3}}{\sqrt{n}}\right\}-1\right).
\end{align}
From the inequality $e^{x}-1 \leq xe^{x}\ (x\geq 0),$
we have
\begin{equation}\label{3.29star}
	\exp\left\{C\frac{1+x^{3}}{\sqrt{n}}\right\}-1\leq C\frac{1+x^{3}}{\sqrt{n}}\exp\left\{C\frac{1+x^{3}}{\sqrt{n}}\right\}.
\end{equation}
Using the inequalities
\begin{equation}\label{3.31}
	\frac{e^{-x^{2} / 2}}{\sqrt{2 \pi}(1+x)} \leq 1-\Phi(x) \leq \frac{e^{-x^{2} / 2}}{\sqrt{\pi}(1+x)}, \quad \ \ x \geq 0,
\end{equation}
and  the inequalities \eqref{3.29} and \eqref{3.29star}, we get
\begin{align}\label{b3.23}
	\mathbb{P}\left(\frac{\log Z_{n}-n \mu}{\sigma \sqrt{n}}\leq x\right)-\Phi(x)
	&\geq	-\frac{C(x^2+1-x)}{\sqrt{2\pi n}}\exp \left\{-\frac{x^{2}}{2}+C\frac{1+x^{3}}{\sqrt{n}}\right\}\notag\\
	&\geq -C\frac{1}{\sqrt{n}}(1+x^{2})\exp \left\{-\frac{x^{2}}{2}+C\frac{x^{3}}{\sqrt{n}}\right\}.
\end{align}
For any $ x\in\left[0,\,n^{1/4}\right] $, it holds
\begin{equation*}
	1-C\frac{x}{\sqrt{n}}
	=\frac{1}{1+\sum\limits_{k=1}^{\infty}(Cx/\sqrt{n})^{k}}
	=\frac{1}{1+Cx/\sqrt{n}\left(\frac{1}{1-Cx/\sqrt{n}}\right)}
	\geq \frac{1}{1+Cx/\sqrt{n}}.
\end{equation*}
Thus, we have for any $ x\in\left[0,\,n^{1/4}\right] $,
\begin{equation}\label{b3.24}
	\exp \left\{-\frac{x^{2}}{2}+C\frac{x^{3}}{\sqrt{n}}\right\}=\exp \left\{-\frac{x^{2}}{2}\left(1-C\frac{x}{\sqrt{n}}\right)\right\}
	\leq\exp\left\{-\frac{x^{2}}{2(1+Cx/\sqrt{n})}\right\}.
\end{equation}
Applying \eqref{b3.24} to \eqref{b3.23}, we obtain for any $ x\in\left[0,\, n^{1/4}\right] $,
\begin{equation}\label{equation10}
	\mathbb{P}\left(\frac{\log Z_{n}-n \mu}{\sigma \sqrt{n}}\leq x\right)-\Phi(x)
	\geq -C\frac{1}{\sqrt{n}}(1+x^{2})\exp\left\{-\frac{x^{2}}{2(1+cx/\sqrt{n})}\right\}.
\end{equation}
Next, we prove the upper bound of $ \mathbb{P}\left(\frac{\log Z_{n}-n \mu}{\sigma \sqrt{n}}\leq x\right)-\Phi(x) $. Again using the first inequality in Lemma \ref{lemma4.9}, we get
\begin{align}\label{a3.38}
	\mathbb{P}\left(\frac{\log Z_{n}-n \mu}{\sigma \sqrt{n}}\leq x\right)-\Phi(x)
	&=-\left[\mathbb{P}\left(\frac{\log Z_{n}-n \mu}{\sigma \sqrt{n}} > x\right)-(1-\Phi(x))\right]\notag\\
	&\leq	-\left[(1-\Phi(x))\exp\left\{-C\frac{1+x^{3}}{\sqrt{n}}\right\}-(1-\Phi(x))\right]\notag\\
	&=-(1-\Phi(x))\left(\exp\left\{-C\frac{1+x^{3}}{\sqrt{n}}\right\}-1\right).
\end{align}
From the inequality $e^{x}-1 \geq x$ $, x\leq 0$,
we have for any $x\leq 0,$
\begin{equation}\label{a3.39}
	\exp\left\{-C\frac{1+x^{3}}{\sqrt{n}}\right\}-1\geq -C\frac{1+x^{3}}{\sqrt{n}}.
\end{equation}
Applying \eqref{3.31} and \eqref{a3.39} to \eqref{a3.38}, we obtain
\begin{align}\label{a3.25}
	\mathbb{P}\left(\frac{\log Z_{n}-n \mu}{\sigma \sqrt{n}}\leq x\right)-\Phi(x)
	&\leq \frac{C(1+x^2-x)}{\sqrt{\pi n}}\exp \left\{-\frac{x^{2}}{2}\right\}\notag\\
	&\leq C\frac{1}{\sqrt{n}}(1+x^{2})\exp\left\{-\frac{x^{2}}{2(1+cx/\sqrt{n})}\right\}.
\end{align}
Combining \eqref{equation10} and \eqref{a3.25} together, we get for any $x\in \left[0,\, n^{1/4}\right]$,
\begin{equation*}
	\left|\mathbb{P}\left(\frac{\log Z_{n}-n \mu}{\sigma \sqrt{n}}\leq x\right)-\Phi(x)\right|
	\leq C\frac{1}{\sqrt{n}}(1+x^{2})\exp\left\{-\frac{x^{2}}{2(1+cx/\sqrt{n})}\right\}.
\end{equation*}

When $ x> n^{1/4}$, we have
\begin{align}\label{equation11}
	\left|\mathbb{P}\left(\frac{\log Z_{n}-n \mu}{\sigma \sqrt{n}}\leq x\right)-\Phi(x)\right| &=\left|\mathbb{P}\left(\frac{\log Z_{n}-n \mu}{\sigma \sqrt{n}} > x\right)-(1-\Phi(x))\right| \notag\\
	&\leq	\mathbb{P}\left(\frac{\log Z_{n}-n \mu}{\sigma \sqrt{n}} > x\right)+1-\Phi(x).
\end{align}
By the first inequality in Lemma \ref{lemma4.8} and Lemma \ref{lemma4.7}, it follows that for any $ x> n^{1/4}$,
\begin{equation}
	\mathbb{P}\left(\frac{\log Z_{n}-n \mu}{\sigma \sqrt{n}} > x\right) \leq C\exp\left\{-\frac{x^{2}}{2(1+cx/\sqrt{n})}\right\}.
\end{equation}
Notice that
\begin{align}
	1-\Phi(x)
	\leq \frac{e^{-x^{2} / 2}}{\sqrt{\pi}(1+x)}\leq \frac{1}{1+x}\exp\left\{-\frac{x^{2}}{2(1+cx/\sqrt{n})}\right\}.
\end{align}
And it holds for any $ x > n^{1/4}$,
\begin{align}\label{a3.43}
	\frac{1}{\sqrt{n}}(1+x^{2})\geq 1\quad\text{and}\quad \frac{1}{1+x}\leq\frac{1}{\sqrt{n}}(1+x^{2}).
\end{align}
Combining \eqref{3.31} and \eqref{equation11}-\eqref{a3.43} together, we get
\begin{equation*}
	\left|\mathbb{P}\left(\frac{\log Z_{n}-n \mu}{\sigma \sqrt{n}}\leq x\right)-\Phi(x)\right|
	\leq C\frac{1}{\sqrt{n}}(1+x^{2})\exp\left\{-\frac{x^{2}}{2(1+cx/\sqrt{n})}\right\},
\end{equation*}
which gives the desired inequality for $ x > n^{1/4}$.

For the case when $ x<0 $, it can be proved in a similar way, but \eqref{3.31} is replaced by
\begin{equation*}
	\frac{e^{-x^{2} / 2}}{\sqrt{2 \pi}(1+|x|)} \leq \Phi(x) \leq \frac{e^{-x^{2} / 2}}{\sqrt{\pi}(1+|x|)},\quad\ x\leq0,
\end{equation*}
and the second inequalities in Lemma \ref{lemma4.8} and  Lemma \ref{lemma4.9} are used for the cases $ x \in\left(-\infty,\, -n^{1/4}\right) $ and $ x\in\left[ -n^{1/4},\,0\right] $, respectively. \hfill$\Box$

\setcounter{equation}{0}

\end{document}